\documentclass[12 pt]{article}
\usepackage[english]{babel}
\usepackage[active]{srcltx}
\usepackage{amsmath}
\usepackage{amsfonts,amssymb}
\usepackage[mathcal]{eucal}
\newtheorem{theorem}{Theorem}

\newtheorem{definition}{Definition}

\newcommand{\rav}{\stackrel{\triangle}{=}}
\newcommand{\ravref}[1]{\stackrel{(\ref{#1})}{=}}
\newcommand{\leqref}[1]{\stackrel{(\ref{#1})}{\leq}}
\newcommand{\geqref}[1]{\stackrel{(\ref{#1})}{\geq}}
\newcommand{\rref}[1]{$(\ref{#1})$}
\newcommand{\mm}[1]{{\mathbb{#1}}}

\newcommand{\ct}[1]{{\mathcal{#1}}}
\newcommand{\epsi}{\varepsilon}
\newcommand{\doc}{{\em{Proof\ }}}
\newcommand{\bo}{\hfill {$\Box$}}
\newcommand{\hjk}[1]{{\tilde{#1}}}

\begin{document}
\title{On transversality condition for overtaking optimality in infinite horizon control problem\thanks{Krasovskii Institute of Mathematics and Mechanics,     Russian
Academy of Sciences, 16, S.Kovalevskaja St., 620990, Yekaterinburg, Russia; \ \
Institute of Mathematics and Computer Science, Ural Federal University, 4, Turgeneva St., 620083, Yekaterinburg, Russia}}

\author{Dmitry Khlopin\\
{\it khlopin@imm.uran.ru} %etc.
}

\maketitle

\begin{abstract}
In this paper we investigate   necessary conditions of optimality for infinite-horizon optimal control problems with overtaking optimality as an optimality criterion. For the case of local Lipschitz continuity of the payoff function, we construct a boundary condition on the co-state arc that is necessary for the optimality. We also show that, under additional assumptions on the payoff function's asymptotic behavior, the Pontryagin Maximum Principle with this condition becomes a complete system of relations,  and this boundary condition points out the unique co-state arc through a Cauchy-type formula. An example is given to clarify the application of this formula as an explicit expression of the co-state arc. The cornerstone of this paper is the theorem on convergence of subdifferentials.


%\subjclass[2010]{}
% 91A25, 49L20, 49N70, 91A23, 40E05

%\keywords{Optimal control \and  Infinite horizon problem \and transversality condition for infinity \and overtaking optimal control \and convergence of subdifferentials
%}
 %\subclass{49J52 \and 49K15 \and  91B62
 %}

\end{abstract}

  \section*{Introduction.}

Necessary conditions on infinite-horizon control problems were proved in their maximally general form by H.Halkin in the form of the Pontryagin Maximum Principle (PMP) in \cite[\S~4]{Halkin}; however, these relations lacked a boundary condition at infinity and could not help to select a unique solution of the adjoint system.

In this paper, we propose a modification of Halkin's general construction of necessary conditions of optimality in which the transversality condition is obtained through the theorems on convergence of subdifferentials. The co-state arc is described through the limiting gradient of the payoff function. For simplicity,  we use the overtaking optimality as the optimality criterion; we also assume the gradients of the payoff function to be bounded (it also implies the normality of the PMP system). We also show that additional conditions imposed on the system---such as the continuous dependence of the payoff function's gradient on the initial conditions---provide for the existence of a unique solution of the PMP system supplemented with the above-mentioned transversality condition. A similar condition was also studied in
     \cite{kr_as}, \cite[\S 4]{kab}, \cite{av_new},  \cite{JDCS}, \cite{Optim}, \cite{Tauchnitz}.
      None of those cover the results of this paper.

%  \section{Base definitions and statements.}\ \\
  \section{Preliminaries.}\ \\

Let
  ${\mathbb{T}}\rav
  {\mathbb{R}}_{\geq 0}$
be the time interval of the initial control system, and let its state space be a certain finite-dimensional Euclidean space ${\mathbb{X}}\rav{\mathbb{R}}^m$.

Consider an infinite-horizon control problem,
\begin{subequations}
\begin{eqnarray}
     l(b)+\int_{0}^\infty f_0(t,x,u)\, dt\to\min  \label{sys0_}\\
    \dot{x}=f(t,x,u),\quad u\in U,  \label{sys_}\\
    x(0)\in\ct{C}.    \label{sysK_}
\end{eqnarray}
\end{subequations}
Here, $r$ and~$f_0$ are scalar functions; $x$ is the state variable, which assumes values from~${\mathbb{X}}$, and~$u$ is some control parameter from a given subset~$U$ of a certain finite-dimensional Euclidean space. Admissible controls are elements of the set ${\ct{U}}\rav L^\infty_{loc}(\mm{T},U).$

We assume the following conditions to hold:
\begin{itemize}
  \item  ${\mathcal{C}}$ is a closed subset of~${\mathbb{X}}$;
  \item   $l$ is a locally Lipschitz continuous scalar function of $x\in\mm{X}$;
  \item for all $u\in {\ct{U}}$, the functions $\mathbb{T}\times\mathbb{X}\ni(t,x)\mapsto f(t,x,u(t))\in\mathbb{X}$ and
   $\mathbb{T}\times\mathbb{X}\ni(t,x)\mapsto f_0(t,x,u(t))\in\mathbb{R}$ and their derivatives with respect to~$x$ are Borel-measurable in~$t$, locally Lipschitz continuous in~$x$, and satisfy the sublinear growth condition with respect to~$x$.
\end{itemize}

Thus, for every admissible control~$u\in {\ct{U}}$, time $\theta\in\mm{T}$,
and initial state $b\in{\mathbb{X}}$, there exists a unique solution $y(b,\theta,u;\cdot)$ of \rref{sys_} with the initial condition $x(\theta)=b$, which can be assumed to be defined for the whole~${\mathbb{T}}$.
Let us now introduce a scalar function~$J$ as follows:
\begin{eqnarray*}
   J(b,\theta;u,T)&\rav&\int_{\theta}^{T} f_0\big(t,y(b,\theta,u;t),u(t)\big)\,dt\quad \forall b\in\mm{X},u\in\ct{U},\theta\in\mm{T},T>\theta.
\end{eqnarray*}
The conditions already imposed guarantee the smoothness of~$J$ in~$x$ and the validity of PMP~\cite[Theorem 5.2.1]{cl_new} for a
finite-horizon control problem.

Call a pair
 $(x,u)\in C(\mm{T},\mm{X})\times {\mathfrak{U}}$ an admissible control process
if $x(0)\in{\mathcal{C}}$ and $x(\cdot)=y(x(0),0,u;\cdot).$
\begin{definition}
Call an admissible process $(\hjk{x},\hjk{u})$ overtaking optimal \cite{car1} for problem \rref{sys0_}--\rref{sysK_} if for every admissible process $(x,u)$ it holds that
\begin{eqnarray*}
%\label{tauover_}
      \liminf_{T\to\infty}\Big(l(x(0))-l(\hjk{x}(0))+
      \int_{0}^T[f_0(t,x(t),u(t))-f_0(t,\hjk{x}(t),\hjk{u}(t))]\,dt
      %\int_0^{T} f_0\big(t,\hjk{x}(t),\hjk{u}(t)\big)dt-
      %-\int_0^{T} f_0\big(t,x(t),u(t)\big)dt
      \Big) \geq 0.
\end{eqnarray*}
\end{definition}
Hereinafter assume that a certain admissible control process $(\hjk{x},\hjk{u})$ is
overtaking optimal for problem \rref{sys0_}--\rref{sysK_}.
For brevity, let us also introduce
$$\hjk{J}(b;T)\rav J(b,0;\hjk{u},T),\quad \hjk{y}(b;T)\rav y(b,0,\hjk{u};T)\qquad\forall T>0,b\in\mm{X}. $$

We will also make use of elementary notions from the nonconvex analysis \cite{lpenot}.
For a Lipschitz continuous function $g:\mm{X}\to \mm{R}$  and a point $\xi\in\mm{X}$,
denote by $\hat{\partial} g(\xi)$  the Fr\'{e}chet subdifferential of this function at the point $\xi\in\mm{X}$;
 it consists of all $h'(\xi)\in\mathbb{X}^*$ for a  function $h : \mm{X}\to {\mathbb{R}}\cup\{+\infty\}$ such that
(a) $h(\xi')\leq g(\xi')$ for every $\xi'\in \mm{X}$ and $g(\xi) = h(\xi)$,
(b) $h$ is Fr\'{e}chet differentiable.
  Denote the limiting subdifferential of $g$ at $\xi$ by $\partial g(\xi)$; it
 consists of all
$\zeta$ in $\mm{X}^*$  such that
\[
\exists \textrm{ sequences of }  y_n\in{\mathbb{X}},\zeta_n\in \hat{\partial} g(y_n),y_n\to \xi,
\zeta_n\to\zeta, g(y_n)\to g(\xi).\]
Denote by $N^{{\mathcal{C}}}(\xi)$  the limiting normal cone of ${\mathcal{C}}$ at $\xi.$

\section{The Pontryagin Maximum Principle and additional transversality conditions}
\label{187}
  Let the Hamilton--Pontryagin function
  $H:{\mathbb{X}}\times{\mathbb{X}}^*\times {U}\times{\mathbb{T}}^2\mapsto{\mathbb{R}}$
  be given by
\begin{eqnarray*}
   H(x,\psi,u,\lambda,t)\rav\psi f\big(t,x,u\big)\!-\!\lambda
   f_0\big(t,x,u\big),\ 
   \forall (x,\psi,u,\lambda,t)\!\in\!{\mathbb{X}}\times{\mathbb{X}}^*\times {U}\times{\mathbb{T}}^2.
\end{eqnarray*}
 Let us introduce the relations of the Pontryagin Maximum Principle:
\begin{subequations}
 \begin{eqnarray}
   \label{sys_x}
       \dot{x}(t)&=& f\big(t,x(t),\hjk{u}(t)\big);\\
   \label{sys_psi}
       -\dot{\psi}(t)&=&\frac{\partial
       {H}}{\partial x}\big(x(t),\psi(t),\hjk{u}(t),\lambda,t\big);\\
   \label{maxH}
\sup_{u'\in
U}H\big(x(t),\psi(t),u',\lambda,t\big)%=\ct{H}\big(x(t),\psi(t),\lambda,t\big)
&=&
        H\big(x(t),\psi(t),\hjk{u}(t),\lambda,t\big).
   \end{eqnarray}
\end{subequations}

From~\cite{Halkin}, it follows that, for an overtaking optimal process, there exists a nontrivial solution of PMP \rref{sys_x}-\rref{maxH}. This system of necessary relations  of optimality lacks one more boundary condition on the adjoint variable, which corresponds to the transversality condition at infinity.

For example, it is possible to construct such a condition if the value function is known, see e.g. \cite{CF},\cite{khlopinPAFA},\cite{sagara}. Another approach is connected with the use of the corresponding Sobolev spaces, see e.g. \cite{aucl},\cite{Tauchnitz}.
The transversality condition that we obtain in this paper is based on the following definition:
\begin{definition}
Call a nontrivial solution $(\hjk{x},\hjk{\psi},\hjk{\lambda})$ of system \rref{sys_x}--\rref{sys_psi} an exact limiting solution iff
for certain sequences of  $y_n\in \mm{X},t_n\in{\mm{T}},\lambda_n>0$
  %$\zeta_n=\frac{\partial \hjk{J}}{\partial x}(y_n;t_n)$
it holds that
     \begin{eqnarray}\nonumber
   t_n\to\infty,y_n\to \hjk{x}(0),\lambda_n\to \hjk{\lambda},\\
\lambda_n\frac{\partial \hjk{J}}{\partial x}(y_n;t_n)\to\hjk{\psi}(0), \hjk{J}(y_n;t_n)-\hjk{J}(\hjk{x}(0);t_n)\to 0. \label{WAK}
     \end{eqnarray}
\end{definition}
As proved in \cite[Proposition~2.1]{Optim}, to every process $(\hjk{x},\hjk{u})$ that is weakly uniformly overtaking optimal \cite{car1} for problem \rref{sys0_}--\rref{sysK_}, one could assign an exact limiting solution $(\hjk{\psi},\lambda^*)$ of PMP \rref{sys_x}--\rref{maxH} with $\hjk{\lambda}\in\{0,1\}$.
See other means of expressing this condition in e.g. \cite{JDCS},\cite{Optim}.

%%YS: "�����������" --> asymptotics? asymptotic behavior? ������� ���������������? ����� ����, "find asymptotic conditions on the adjoint system that would hold for at least a single solution but would not hold for a continuum of solution? ��� ������� that would yield at least one but at most a countable number of solutions

In infinite-horizon control problems, a principal obstacle to obtaining additional conditions, the transversality conditions, is the need to
find asymptotic conditions on the adjoint system that would hold for at least a single solution but would not hold for a continuum of solutions.
%find such an expression of asymptotic behavior of the adjoint system that would yield at least one and at most a countable number of solutions.
 In certain problems, it is possible to find a condition that assigns to each optimal process exactly one solution of the adjoint system.
 To spell the formula that describes this condition, let us first recall the Cauchy formula for adjoint systems.

%%YS: �������� ������� \operatoraname{dim}; �������� ��������� �������� ����� ��� ������������ ����������; ������, �����, ������ ����� ������
% OK
 Denote by~${\mathbb{L}}$ the linear space of all real $m\times m$ matrices; here, $m=\operatorname{dim}{\mathbb{X}}$.
For each $\xi\in{\mathbb{X}}$, there exists a solution ${A}(\xi;t)\in C({\mathbb{T}}, {\mathbb{L}})$ of the Cauchy problem
 \begin{equation*}
 \frac{d{A}(\xi;t)}{dt} =\frac{\partial f }{\partial x}
 \big(\hjk{y}(\xi;t),\hjk{u}(t)\big)
  A(\xi;t),\quad A(\xi;0)=1_{\mathbb{L}}.
\end{equation*}
Then,
  \begin{equation}
   \label{4C}
\frac{\partial \hjk{y}}{\partial x}(\xi;T)=A(\xi;t),\quad
\frac{\partial \hjk{J}}{\partial x}(\xi;T)%=\frac{\partial J}{\partial x}(x(0),0;\hjk{u},T)
%=I(x(0),T).
=\int_0^T
   \frac{\partial f_0}{\partial x}
    \big(t,\hjk{y}(\xi;t),\hjk{u}(t)\big)
\, A(\xi;t)
  \,dt
\end{equation}
 and, for every~$\lambda$, its   solution $(x,\psi)$ of system \rref{sys_x}--\rref{sys_psi}
 satisfies  the following Cauchy formula:
  \begin{equation}
   \label{4A}
   \psi(t)A(x(0);t)-\psi(0)=
   \lambda \frac{\partial \hjk{J}}{\partial x}(x(0);T)
%   \lambda I(x(0);t)
   \qquad \forall t\in{\mathbb{T}}.
\end{equation}

In papers
     \cite{kr_as,kab}, and then in \cite{av_new,belyakov,Tauchnitz}, a number of assumptions on the asymptotic behavior of $f,f_0,J$, and their derivatives was obtained, which provide for a unique reconstruction of the PMP solution (through $(\hjk{x},\hjk{u})$) by means of the formulas
  \begin{equation} \label{AK}
  -\hjk{\psi}(0)=\lim_{T\to\infty} \frac{\partial \hjk{J}}{\partial x}(\hjk{x}(0);T)=\int_0^\infty
   \frac{\partial f_0}{\partial x}
    \big(t,\hjk{x}(t),\hjk{u}(t)\big)
\, A(\hjk{x}(0);t)
  \,dt,\  \lambda^*=1.
\end{equation}
We study the possibility of using conditions   \rref{WAK} and~\rref{AK} assuming only the boundedness of $\frac{\partial \hjk{J}}{\partial x}$.
	In addition, based on condition~\rref{AK}, we will also prove the necessity of another, supplementary condition:
for all $u\in U$ and almost all $t\geq 0,$
\begin{eqnarray}\nonumber
   \liminf_{T\to\infty}\Big[{H}\big(\hjk{x}(t),-\frac{\partial J}{\partial x}(\hjk{x}(t),t;\hjk{u},T),\hjk{u}(t),1,t\big)&\ &\\
   -
   {H}\big(\hjk{x}(t),-\frac{\partial J}{\partial x}(\hjk{x}(t),t;\hjk{u},T),u,1,t\big)\Big]&\geq& 0. \label{Anton}
\end{eqnarray}
%%YS: ����� an ... control, ������ ��� � �� ����, ���������� �� �� ������������; ���� ����������, �� the ... control
 Such a condition was proposed in~\cite{belyakov} as a means of seeking an overtaking optimal control.

\section{The main result}

%\begin{subequations}
\begin{theorem}
\label{4}
    Let the process $(\hjk{x},\hjk{u})$ be overtaking  optimal for \rref{sys0_}--\rref{sysK_}.

 Assume that, for every bounded neighborhood~$\Xi$  of the point $\hjk{x}(0)$, for all $T>0,\xi\in\Xi$,
the vectors
 $\frac{\partial \hjk{J}}{\partial x}(\xi;T)=\frac{\partial J}{\partial x}(\xi,0;\hjk{u},T)$
are uniformly bounded.

   Then,  there exists an exact limiting solution $(\hjk{x},\hjk{\psi},1)$ of PMP \rref{sys_x}-\rref{maxH}
such that
\begin{eqnarray}
\label{400}
%\textrm{(transversality  condition at zero)} \qquad\qquad
 \hjk{\psi}(0)\in \partial l(\hjk{x}(0))+N^{{\mathcal{C}}}(\hjk{x}(0)),
\end{eqnarray}
 in particular, $-\hjk{\psi}(0)$
 is a partial limit of $\frac{\partial J}{\partial x}(\xi,0;\hjk{u},T)$
  {as} $\xi\to \hjk{x}(0),T\to\infty.$
\end{theorem}
\begin{theorem}\label{5}
Under conditions of Theorem~\ref{4}, let there also exist a finite limit
\begin{eqnarray}
%\label{400}
\label{lim}
    \lim_{\xi\to\hjk{x}(0),T\to\infty}\frac{\partial J}{\partial x}(\xi,0;\hjk{u},T).
\end{eqnarray}
%   b) ������� $J(\cdot,0;\hjk{u},T)$ ��������� �������� ��� $T\to\infty$ �
%   ��� ������� $\epsi>0$  � ��������� ����������� $\Xi_\epsi$ ����� $\hjk{x}(0)$ ��� ���� $T>0$,$\xi,\xi'\in\Xi_\epsi$ ���������
%   $$J(\xi',0;\hjk{u},T)-J(\xi,0;\hjk{u},T)-\frac{\partial J}{\partial x}(\xi,0;\hjk{u},T)(\xi-\xi')\leq \epsi||\xi-\xi'||.$$

%%YS: "converges in the Riemann sense" --- ����� 2 ��������� �� Scholar

	Then, %!the limit in~\rref{AK} exists and finite!, and
  the system of relations \rref{sys_x}-\rref{maxH},\rref{AK} has exactly one solution. Moreover, this solution also satisfies condition~\rref{Anton}.
\end{theorem}

These propositions are all proved in the next section.

Let us show that if condition~\rref{lim} does not hold, then, under conditions of Theorem~\ref{4}, formula \rref{AK} may not specify a solution of PMP \rref{sys_x}-\rref{maxH}. To this end, consider an example where all the maps $x\mapsto \hjk{J}(x,T)$ are 1-Lipschitz continuous, however, condition~\rref{AK} specifies the solution of system \rref{sys_x}-\rref{sys_psi} that does not satisfy condition of maximality of the Hamiltonian~\rref{maxH}.

Consider the following problem:
%\begin{subeq1uations}
\begin{eqnarray*}
     \int_{1}^2 \frac{1}{2}\sin (2x)\, dt+\int_{2}^\infty\Big[\frac{x}{t}\cos (tx)-\frac{1}{t^2}\sin (tx)\Big]\, dt\to\min\\
    \dot{x}=u 1_{[0,1]}(t),\quad u\in [-1,1],\\
    x(0)=0.
\end{eqnarray*}
%\end{subequations}
    Let us first look at the map $M(T,x)\rav\frac{1}{T}\sin (Tx)$ for $T\geq 1,x\in \mm{R}$; it is differentiable, and its partial derivatives equal, respectively,
     $$\frac{\partial M}{\partial T}(T,x)=\frac{x}{T}\cos (Tx)-\frac{1}{T^2}\sin (Tx),\qquad \frac{\partial M}{\partial x}(T,x)=\cos (Tx);$$
in particular, it is 1-Lipschitz continuous in~$x.$

      Consider an arbitrary admissible process $({x},{u})$. For it, we have $x|_{[1,\infty)}\equiv x(1).$
      Now,
      \begin{eqnarray*}
J(0,0;u,2)=\frac{1}{2}\sin (2x(1))=M(2,x(1))=M(2,x(2)),\\
      J(0,0;u,T)=M(2,x(1))+M(T,x(T))-M(2,x(2))=M(T,x(T)).
\end{eqnarray*}
      Since $|J(0,0;u,T)|\leq\frac{1}{T}$ for $T>2$, every admissible process $({x},{u})$ is overtaking optimal (moreover, strongly optimal \cite{car1} and classical optimal \cite{slovak}).
	  Thus, for all admissible processes $(\hjk{x},\hjk{u})$, all conditions of Theorem~\ref{4} hold.

	  Let us prove that, for the overtaking optimal process $(\hjk{x},\hjk{u})\equiv 0$, the implication of Theorem~\ref{5} does not hold.
      Clearly, for all $\xi\in\mm{R}, T>1$, we have
       $$J(\xi,0;\hjk{u},T)=M(T,\xi),\quad \frac{\partial J}{\partial x}(\xi,0;\hjk{u},T)=\cos (T\xi),\quad
       \frac{\partial J}{\partial x}(\hjk{x}(0),0;\hjk{u},T)=1.$$
			By Theorem~\ref{5}, there should exist a solution $(\hjk{x},\hjk{\psi},1)$ of relations \rref{sys_x}-\rref{maxH} that satisfies the initial condition $\hjk{\psi}(0)=-1$.

       Since $H(x,\psi,u,\lambda,t)\rav\psi u$ for all $x,\psi,u,\lambda\in\mm{R},t\in(0,1)$, by \rref{sys_psi}, this fact would imply that $\hjk{\psi}|_{[0,1]}\equiv -1$, and, from \rref{maxH}, $$0=H(\hjk{x}(t),\hjk{\psi}(t),\hjk{u}(t),1,t)=\max_{u\in[-1,1]}H(\hjk{x}(t),\hjk{\psi}(t),u,1,t)=1.$$
 		The obtained contradiction proves that, in the considered example for $(\hjk{x},\hjk{u})\equiv 0$, the result of Theorem~\ref{5} does not hold; therefore, condition \rref{lim} can not be excluded from the conditions of Theorem~\ref{5}.

\section{Theorem proofs
 }
%%YS: ������������� ���������� --> uniformly equicontinuous?
%%YS: ������������� --> precompact? �������, ��� ����������� ������������� ��-�������, �����, ����������������

\doc of Theorem~\ref{4}.
Since, for every bounded neighborhood~$\Xi$ of the point~$\hjk{x}(0)$, the mappings
  $$\Xi\ni\xi\mapsto \frac{\partial \hjk{J}}{\partial x}(\xi;T),\qquad \Xi\ni\xi\mapsto \frac{\partial \hjk{J}}{\partial x}(\xi;T)-\frac{\partial \hjk{J}}{\partial x}(\hjk{x}(0);T)\qquad \forall T>0$$
are uniformly (in $T>0$) bounded, the mappings
  $$\Xi\ni\xi\mapsto \hjk{J}(\xi;T)-\hjk{J}(\hjk{x}(0);T)\qquad \forall T>0$$
share a common Lipschitz constant~$L$; they are also uniformly equicontinuous.
Since all these mappings become zero at $\xi=\hjk{x}(0)$, they are also uniformly bounded, therefore, the family of these mappings is precompact. Hence, the closure of
  $\{\mm{X}\ni\xi\mapsto \hjk{J}(\xi;T)-\hjk{J}(\hjk{x}(0);T)\,|\, T>0\}$
%  !is embedded into a compact subset of $C(\mm{X},\mm{R})$!
   %is also pre
   is  compact
in the compact-open topology.

Fix an arbitrary unboundedly increasing sequence of positive~$t_n$.
Removing some elements if necessary, it is safe to assume that the mappings $\mm{X}\ni\xi\mapsto \hjk{J}(\xi;t_n)-\hjk{J}(\hjk{x}(0);t_n)$
converge to a certain locally Lipschitz continuous mapping uniformly on every compact. Note that for all $\xi_1\in\mm{X}$ and $t>0$ there exists  $\xi\in\mm{X}$ such that $\xi_1=\hjk{y}(\xi,t)$.
Since for all $\xi\in\mm{X}$, $t\in\mm{T}$, and $T>t$ we have
\begin{eqnarray}\label{420}
 \hjk{J}(\xi,T)-\hjk{J}(\xi,t)=J(\hjk{y}(\xi,t),t;\hjk{u},T),
\end{eqnarray}
there exists the following limit:
\begin{eqnarray}\label{nado}
 {J}_*(\xi,t)&\rav&
% \inf_{k\in\mm{N}}
 \lim_{k\to\infty} \big[J(\xi,t;\hjk{u},\tau_k)-J(\hjk{x}(t),t;\hjk{u},\tau_k)\big],\qquad \forall t\geq 0, \xi\in\mm{X},
\end{eqnarray}
which is uniform in every compact subset of the set~$\mm{X}$.
Note that the mapping~${J}_*$ is also locally Lipschitz continuous; moreover, for all $\xi\in\mm{X}$ and $T>0$, we have
\begin{eqnarray}\nonumber
 {J}_*(\xi,0)-{J}_*(\hjk{y}(\xi,T),T)&=&\lim_{k\to\infty}\big[
 J(\xi,0;\hjk{u},\tau_k)-J(\hjk{x}(0),0;\hjk{u},\tau_k)\\
 &\ &\qquad\quad-J(\hjk{y}(\xi,T),T;\hjk{u},\tau_k)+J(\hjk{x}(T),T;\hjk{u},\tau_k)\big]\nonumber\\
 &\ravref{420}&\hjk{J}(\xi,T)-\hjk{J}(\hjk{x}(0),T).\label{422}
\end{eqnarray}
%%YS: ������ ��� ����� ���������? ���������� ������� �����������������
%% Thank you!!!
For arbitrary $T>0$, denote by  $\hat{\partial}_x J_*(\cdot,T)$, $ \partial_x J_*(\cdot,T)$ the corresponding subdifferentials of the mappings  $\xi\mapsto {J}_*(\xi;T)$.

Since the mapping $\xi\to y(\xi;t)$ is a diffeomorphism for arbitrary~$t$
and the mapping $\xi \to \hjk{J}(\xi;t)$ is smooth, from the elementary properties of the limiting subdifferential (see e.g. \cite[Proposition 6.17]{lpenot}), it follows that
%     Penot, Jean-Paul. Calculus without derivatives. Vol. 266. Springer Science & Business Media, 2012.
 \begin{eqnarray}
  \frac{\partial \hjk{J}}{\partial \xi}(\xi;T)
  &\ravref{422}&\frac{\partial }{\partial \xi} \Big(J_*(\xi;0)-J_*(\hjk{y}(\xi,T),T)+\hjk{J}(\hjk{x}(0),T)\Big),\nonumber\\
  &=&\frac{\partial }{\partial \xi} \Big(J_*(\xi;0)-J_*(\hjk{y}(\xi,T),T)\Big),\nonumber\\
%  \hat{\partial}_x J_*(\xi,0)&\ravref{4C}&\frac{\partial \hjk{J}}{\partial \xi}(\xi;T)+
%  \hat{\partial}_x J_*(\hjk{y}(\xi,T),T)A(\xi;T),\nonumber\\
  \partial_x J_*(\xi,0)&=&\frac{\partial \hjk{J}}{\partial \xi}(\xi;T)+
  \partial_x J_*(\hjk{y}(\xi,T),T)A(\xi;T)\quad \forall\xi\in\mm{X},T>0. \label{424}
  \end{eqnarray}

By virtue of overtaking optimality of~$\hjk{u}$, we have
    $$\liminf_{n\to\infty} \Big[l(b)+J(b,0;u,t_n)-\hjk{J}(\hjk{x}(0);t_n)\Big]\geq l(\hjk{x}(0))\qquad\forall u\in\ct{U}, b\in\ct{C}.$$
Then the same also holds true for  $u\in\ct{U}$ such that $u|_{[t_n,\infty)}=\hjk{u}|_{[t_n,\infty)}$ for a certain $n\in\mm{N}$, whence
\begin{eqnarray*}
   l(\hjk{x}(0))&\leq&\liminf_{k\to\infty} \Big[ l(b)+J(b,0;u,t_n)+J(\hjk{y}(b;t_n),t_n;\hjk{u},t_k)-\hjk{J}(\hjk{x}(0);t_k)\Big]\\
     %&=&l(b)+J(b,0;u,t_k)-\hjk{J}(\hjk{x}(0);t_k)+\liminf_{n\to\infty} J(\hjk{y}(b;t_k),t_k;\hjk{u},t_n)-J(\hjk{x}(t_k),t_k;\hjk{u},t_n)\\
     &=&l(b)+J(b,0;u,t_n)-\hjk{J}(\hjk{x}(0);t_n)+{J}_*(\hjk{y}(b;t_n),t_n).
\end{eqnarray*}
for all  $u\in\ct{U}$, $b\in\ct{C}$, and $n\in\mm{N}.$

%%YS: ����������� �������� --> optimal value

Therefore, for every  $n\in\mm{N}$,    the optimal value of the problem
\begin{eqnarray*}
     l(x(0))+\int_0^{t_n} \Big[f_0\big(t,x(t),u(t)\big)-f_0\big(t,\hjk{x}(t),\hjk{u}(t)\big)\Big]dt
       +{J}_*(x(t_n),t_n)\to\min  %\label{sys0_1}
       \\
    \dot{x}=f(t,x,u),\quad u\in U,  %\label{sys_1}
    \\
    x(0)\in\ct{C}    %\label{sysK_1}
\end{eqnarray*}
is not less than $l(\hjk{x}(0))$.
Consequently, $(\hjk{x},\hjk{u})$ is optimal in such a problem for arbitrary natural~$n$.

Now, for every $n\in\mm{N}$, by \cite[Theorem~5.2.1]{cl_new}, there exist
 ${\psi}_n\in C({\mathbb{T }},{\mathbb{X}}^*)$ such that every triple  $(\hjk{x},\psi_n,1)$ satisfies PMP \rref{sys_x}-\rref{maxH}
    almost everywhere in $[0,t_n]$ with the boundary conditions
\begin{eqnarray}
\psi_n(0)&\in& \partial l (\hjk{x}(0))+N^{\ct{C}}(\hjk{x}(0)),\label{trans_0_max_}\\
%\psi_n(0)&\in& \partial^1 l(\td{x}_n(0))+N_L^{{\mathcal{C}}}(\hjk{x}(0)),
%\label{trans_0_max_}\\
-\psi_n(t_n)&\in& \partial_x J_*(\hjk{x}(t_n),t_n).\label{trans_n_max_}
   \end{eqnarray}

In particular, $\psi_n$, as a solution of \rref{sys_psi}, satisfies  the Cauchy formula (see \rref{4A}), and, by a sequential application of
\rref{4A},\rref{trans_n_max_}, and \rref{424}, we obtain
 \begin{eqnarray*}
 -\psi_n(0)&\ravref{4A}&
-\psi_n(t_n)A(\hjk{x}(0);t_n)+\frac{\partial \hjk{J}}{\partial \xi}(\hjk{x}(0);t_n)\\
&\in& \partial_x J_*(\hjk{x}(t_n),t_n)A(\hjk{x}(0);t_n)+\frac{\partial \hjk{J}}{\partial \xi}(\hjk{x}(0);t_n)\\
&\ravref{424}&\partial_x J_*(\hjk{x}(0),0)-\frac{\partial \hjk{J}}{\partial \xi}(\hjk{x}(0);t_n)+
\frac{\partial \hjk{J}}{\partial \xi}(\hjk{x}(0);t_n)=\partial_x J_*(\hjk{x}(0),0),%\label{trans_n_max_}
 \end{eqnarray*}
therefore,  $-\psi_n(0)\in\partial_x J_*(\hjk{x}(0),0).$

Since~$J_*$ is locally Lipschitz continuous in~$x$, we have proved the boundedness of the vectors~$\psi_n(0).$
Passing from the sequence of~$t_n$ to its certain subsequence if necessary, we can assume that the sequence of $\psi_n(0)$ converges. Hence, by the theorem on continuous dependence of differential equations' solutions on initial conditions, the sequence of~$\psi_n$ converges in  $[0,\infty)$ to a certain solution~$\hjk{\psi}$ of \rref{sys_psi},
and this convergence is uniform in arbitrary compact time intervals. But, consequently, the triple $(\hjk{x},\hjk{\psi},1)$ also  satisfies relations \rref{sys_x}-\rref{maxH} on the whole $\mm{T}$; moreover, now,  for~$\hjk{\psi}$, condition~\rref{400} is implied by \rref{trans_0_max_}, and $-\psi_n(0)\in\partial_x J_*(\hjk{x}(0),0)$ yields
$-\hjk{\psi}(0)\in\partial_x J_*(\hjk{x}(0),0).$

It remains to prove that $(\hjk{x},\hjk{\psi},1)$  is an exact limiting solution of \rref{sys_x}-\rref{sys_psi}.
Recall that~$J_*$ is a limit of the sequence of mappings $\xi\mapsto \hjk{J}(\xi;t_n)$
that is uniform in a certain neighborhood of the point $\hjk{x}(0)$.
As showed in \cite[Theorem 6.1]{ledtrie}, this means that every element from the Fr\'{e}chet subdifferential $\hat{\partial}_x J_*(z,0)$ (for all $z\in\mm{X}$) can be rendered as a limit of
  $\frac{\partial \hjk{J}}{\partial x}(\xi_i;t_{n(i)})=\frac{\partial J}{\partial x}(\xi_i,0;\hjk{u},t_{n(i)})$ for certain sequences  $\xi_i\to z,n(i)\to\infty$.
 By the definition of the limiting subdifferential, every element from  ${\partial}_x J_*(z,0)$ (for all $z\in\mm{X}$) can be expressed---in view of a certain converging to $z$ sequence of $\xi_i$---as a limit of elements from $\hat{\partial}_x J_*(\xi_i,0)$; however, it implies that every element of ${\partial}_x J_*(z,0)$ is a limit of
  $\frac{\partial \hjk{J}}{\partial x}(\xi_i;t_{n(i)})=\frac{\partial J}{\partial x}(\xi_i,0;\hjk{u},t_{n(i)})$
for certain subsequences of $\xi_i\to z$, $n(i)\to\infty.$
By $-\hjk{\psi}(0)\in\partial_x J_*(\hjk{x}(0),0)$,
there exist a sequence of~$\xi_i$ that converges to $\hjk{x}(0)$
and an unboundedly increasing sequence of natural $n(i)$ such that
     $-{\psi}^*(0)=\lim_{i\to\infty}\frac{\partial J}{\partial x}(\xi_{i},0;t_{n(i)}).$
Since the mappings  $\Xi\ni\xi\mapsto\hjk{J}(\xi;t)$ all have equal Lipschitz constants in an arbitrary bounded domain $\Xi\subset \mm{X}$, from
     $||\xi_i-\hjk{x}(0)||\to0$, it automatically follows that $|J(\xi_{i},0;t_{n(i)})-J(\hjk{x}(0),0;t_{n(i)})|\to 0$. Thus, the triple  $(\hjk{x},\hjk{\psi},1)$  is an exact limiting solution of PMP, which is what we wanted to prove.
   \bo

   \doc of Theorem~\ref{5}.

In~\rref{AK}, the existence and finiteness of the integral  is an immediate consequence of~\rref{4C}. By means of Theorem~\ref{4},  we can pick a solution $(\hjk{x},\hjk{\psi},1)$ of PMP \rref{sys_x}-\rref{maxH} such that  $-\hjk{\psi}(0)$ is a partial limit of $\frac{\partial \hjk{J}}{\partial x}(\xi_n;t_n)$
for certain sequences $\xi_n\to\hjk{x}(0), t_n\to\infty.$
Then, by \rref{lim}, it is also a limit of $\frac{\partial \hjk{J}}{\partial x}(\hjk{x}(0);t)$ as $t\to\infty.$ Now, from
 \rref{4C}, we see that \rref{AK} holds for $(\hjk{x},\hjk{\psi},1)$.
Note that condition~\rref{AK} lets us reconstruct $(\hjk{x},\hjk{\psi},1)$ uniquely.
At the same time, \rref{maxH} holds for all $t\geq 0$ except a possibly empty subset $N\subset\mm{R}$ of measure zero. Fix this set.

Let us prove condition~\rref{Anton}. Suppose it is false. Then, for a certain $\tau\in \mm{R}\setminus N$ and a certain $u\in P$ there exist an unboundedly increasing sequence of times $t'_n$ and a positive number~$\epsi$ such that
\begin{eqnarray}\label{antin}
   {H}\Big(\hjk{x}(\tau),-\frac{\partial J}{\partial x}(\hjk{x}(\tau),\tau;\hjk{u},t'_n),\hjk{u}(\tau),1,\tau\Big)\leq\\
   {H}\Big(\hjk{x}(\tau),-\frac{\partial J}{\partial x}(\hjk{x}(\tau),\tau;\hjk{u},t'_n),u,         1,\tau\Big)-\epsi,\nonumber
\end{eqnarray}
Passing to a subsequence, we can again assume the sequence of vectors  $\frac{\partial J}{\partial x}(\hjk{x}(0),0;t'_n,\hjk{u})$ to converge as $n\to\infty$ and the mappings $\xi\mapsto J(0,\xi;\hjk{u},t'_n)$ to converge uniformly in arbitrary compacts. Repeating the reasoning above for the sequence of $t_n=t'_n$, we obtain $\hjk{\psi}(\cdot)$ as a pointwise limit of $\frac{\partial J}{\partial x}(\hjk{x}(\cdot),\cdot;\hjk{u},t'_n)$.
Passing to the limit in \rref{antin}, for $\tau\in \mm{R}\setminus N$, we have
\begin{eqnarray*}
    {H}\big(\hjk{x}(\tau),\hjk{u}(\tau),\hjk{\psi}(\tau),1,\tau\big)
    \leq
   {H}\big(\hjk{x}(\tau),u,\hjk{\psi}(\tau),1,\tau\big)-\epsi,
\end{eqnarray*}
which contradicts condition~\rref{maxH}, whereas \rref{maxH} holds for the whole $\mm{R}\setminus N$. Condition \rref{Anton} is proved.
    \bo

  \section*{Acknowledgements}
%   I would like to
%express my gratitude to  Ya.V. Salii for
%the translation.
 I would like to express my gratitude
to B.Mordukhovich for valuable discussion in
course of writing this article.
This study was partially supported by the Russian Foundation for Basic Research, project no.
16-01-00505.

\end{document}